\documentclass[3p,times]{elsarticle}

\usepackage{ecrc}

\volume{00}

\firstpage{1}

\journalname{Applied Mathematics and Computation}

\runauth{U.N. Katugampola}

\jid{amc}

\jnltitlelogo{\small{Applied Mathematics\\ and\\ Computation}}

\usepackage{amssymb}
\usepackage{amsthm}
\usepackage{epstopdf}
\usepackage{verbatim}
\usepackage{amsmath,amssymb}
\usepackage{slashbox}
\usepackage{color,soul} 
\usepackage[normalem]{ulem}  
\usepackage{textcomp}   
\usepackage{lscape}    
\usepackage{hyperref}
\usepackage{subfigure}
\newcommand{\dx}[1]{\ensuremath{x^{#1}\frac{d}{dx}}}
\newcommand{\xdx}[1]{x\frac{d^{#1}}{dx^{#1}}}
\newcommand{\ddx}[2]{\ensuremath{x^{#1}\frac{d^{#2}}{dx^{#2}}}}

\newcommand{\ipa}{\ensuremath{[Re(\alpha)]+1}}

\numberwithin{equation}{section}

\newtheorem{theorem}{Theorem}[section]

\newtheorem{lemma}[theorem]{Lemma}
\newtheorem{remark}{Remark}
\newtheorem{definition}[theorem]{Definition}

\newtheorem{algorithm}{Algorithm}

\biboptions{sort&compress}

\usepackage[figuresright]{rotating}

\begin{document}

\begin{frontmatter}


\tnotetext[label1]{\emph{Email~address:} \mbox{uditanalin@yahoo.com}  \\ \mbox{}\hspace{.4cm} Submitted on Jan 28, 2014  \hspace{0.7cm} Revised on Oct 28, 2014  \hfill Accepted for publication in Applied Mathematics and Computation}


\title{Mellin Transforms of Generalized Fractional Integrals and Derivatives}


\author{Udita N. Katugampola}

\address{Department of Mathematical Sciences, University of Delaware, Newark DE 19716, USA.}

\begin{abstract}
We obtain the Mellin transforms of the generalized fractional integrals and derivatives that generalize the Riemann-Liouville and the Hadamard fractional integrals and derivatives. We also obtain interesting results, which combine generalized $\delta_{r,m}$ operators with generalized Stirling numbers and Lah numbers. For example, we show that $\delta_{1,1}$ corresponds to the Stirling numbers of the $2^{nd}$ kind and $\delta_{2,1}$ corresponds to the unsigned Lah numbers. Further, we show that the two operators $\delta_{r,m}$ and $\delta_{m,r}$, $r,m\in\mathbb{N}$, generate the same sequence given by the recurrence relation
\[
S(n,k)=\sum_{i=0}^r \big(m+(m-r)(n-2)+k-i-1\big)_{r-i}\binom{r}{i} S(n-1,k-i), \;\;  0< k\leq n,
\]
with $S(0,0)=1$ and $S(n,0)=S(n,k)=0$ for $n>0$ and $1+min\{r,m\}(n-1) < k $ or $k\leq 0$. Finally, we define a new class of sequences for $r \in \{\frac{1}{3}, \frac{1}{4}, \frac{1}{5}, \frac{1}{6}, \ldots \}$ and in turn show that $\delta_{\frac{1}{2},1}$ corresponds to the generalized Laguerre polynomials. 
\end{abstract}


\end{frontmatter}



\section{Introduction}
The Fractional Calculus (FC) is the generalization of the classical calculus to arbitrary orders. The history of FC goes back to seventeenth century, when in 1695 the derivative of order $\alpha = \frac{1}{2}$ was described by Leibniz. Since then, the new theory turned out to be very attractive to mathematicians as well as biologists, economists, engineers and physicists \cite{what}. In \cite{key-9}, Samko et al. provide a comprehensive treatment of the subject. Several different derivatives were studied: Riemann-Liouville, Caputo, Hadamard, Erd\'{e}lyi-Kober, Grunwald-Letnikov and Riesz are just a few to name \cite{key-8,key-9}. In \cite{udita2}, the author provides an extended reference to some of these and other fractional derivatives.  

In fractional calculus, the fractional derivatives are defined via fractional integrals \cite{key-9, key-11, key-8}. According to the literature, the Riemann-Liouville fractional derivative, hence the Riemann-Liouville fractional integral has played a major roles in FC \cite{key-9}. The Caputo fractional derivative has also been defined via the Riemann-Liouville fractional integral \cite{key-9}. Another important fractional operator is the Hadamard operator \cite{key-10,key-2}. In \cite{key-6}, Butzer, et al. obtained the Mellin transforms of the Hadamard integral and differential operators. The interested reader is referred to, for example, \cite{key-1,key-5,key-6,key-3,key-7,key-8,key-2,key-9} for further properties of those operators. 

In \cite{key-0}, the author introduces a new fractional integral, which generalizes the Riemann-Liouville and the Hadamard integrals into a single form. The interested reader is referred, for example, to \cite{Herr, u-1,u-2,u-3,u-4,u-5,u-6,u-7,u-8,u-9,u-10} for further results on these and similar operators. In \cite{udita2}, the author introduces a new fractional derivative which generalizes the the Riemann-Liouville and the Hadamard fractional derivatives to a single form. 
The present work is primarily devoted to the Mellin Transforms of the generalized fractional operators developed in \cite{key-0} and \cite{udita2}. 

The paper is organized as follows. In the next section, we give definitions and some basic properties of the fractional integrals and derivatives of various types. In section 3, we develop the Mellin transforms of the generalized fractional integrals and derivatives. Further, we investigate the relationship between the Mellin transform operator and the generalized Stirling numbers of the $2^{nd}$ kind. Finally, we introduce a new class of sequences. 

\section{Definitions}
The \emph{Riemann-Liouville fractional integrals} $I^{\alpha}_{a+}f$ and $I^{\alpha}_{b-}f$ of order $\alpha \in \mathbb{C},$ $(Re(\alpha)>0)$ are defined by \cite{key-9},
\begin{equation}
(I^\alpha_{a+} f)(x) = \frac{1}{\Gamma(\alpha)}\int_a^x (x - \tau)^{\alpha -1} f(\tau) d\tau  \quad ; x > a.
\label{eq:lRLI}
\end{equation}
and 
\begin{equation}
(I^\alpha_{b-} f)(x) = \frac{1}{\Gamma(\alpha)}\int_x^b (\tau - x)^{\alpha -1} f(\tau) d\tau  \quad ; x < b,
\label{eq:rRLI}
\end{equation}
respectively, where $\Gamma(\cdot)$ is the Gamma function. These integrals are called the \emph{left-sided} and \emph{right-sided fractional integrals}, respectively. When $\alpha = n \in \mathbb{N},$ the integrals (\ref{eq:lRLI}) and (\ref{eq:rRLI}) coincide with the $n$-fold integrals \cite[chap.2]{key-8}. 
The \emph{Riemann-Liouville fractional derivatives} (RLFD) $D^{\alpha}_{a+}f$ and $D^{\alpha}_{b-}f$ of order $\alpha \in \mathbb{C},\, (Re(\alpha)\geq 0)$ are defined by\cite{key-9},
\begin{equation}
(D^\alpha_{a+} f)(x)=\bigg(\frac{d}{dx}\bigg)^n \, \Big(I^{n-\alpha}_{a+} f\Big)(x) \quad x>a,
\label{eq:lRLD}
\end{equation}
and
\begin{equation}
\;\;\;\;(D^\alpha_{b-} f)(x)=\bigg(-\frac{d}{dx}\bigg)^n \, \Big(I^{n-\alpha}_{b-} f\Big)(x)  \quad x<b,
\label{eq:rRLD}
\end{equation}
respectively, where $n = \ipa$, which is sometimes denoted also by the ceiling function, $\lceil Re(\alpha) \rceil$, when there is no room for confusion. Here $[\cdot]$ represents the integer part. For simplicity, from this point onwards, except in few occasions, we consider only the \emph{left-sided} integrals and derivatives. The interested reader may find more detailed information about the \emph{right-sided} integrals and derivatives in the references, for example in \cite{key-8,key-9}. 

The next type that we elaborate in this paper is the \emph{Hadamard Fractional integral} \cite{key-10,key-8,key-9} given by,
\begin{equation}
\mathbf{I}^\alpha_{a+} f(x) = \frac{1}{\Gamma(\alpha)}\int_a^x \Bigg(\log\frac{x}{\tau}\Bigg)^{\alpha -1} f(\tau)\frac{d\tau}{\tau} \quad ;Re(\alpha) > 0, \, x > a\geq 0.
\label{eq:HI}
\end{equation}
The corresponding \emph{Hadamard fractional derivative} of order $\alpha \in \mathbb{C}, Re(\alpha) >0$ is given by,
\begin{equation}
\mathbf{D}^\alpha_{a+} f(x) = \frac{1}{\Gamma(n-\alpha)}\Bigg(x\frac{d}{dx}\Bigg)^n \int_a^x \Bigg(\log\frac{x}{\tau}\Bigg)^{n-\alpha +1} f(\tau)\frac{d\tau}{\tau} \quad ;x > a\geq 0.
\label{eq:HD}
\end{equation}
where $n=\lceil Re(\alpha) \rceil$. 

In 1940, important generalizations of the Riemann-Liouville fractional operators were introduced by Erd\'{e}lyi and Kober \cite{E-K}. \emph{Erd\'{e}lyi-Kober}-type fractional integral and differential operators \cite{E-K,key-4,key-8,key-13,key-9,key-Y} are defined by,

\begin{align}
\;\;(I^\alpha_{a+; \,\rho, \,\eta} f)(x) &= \frac{\rho\,x^{-\rho(\alpha+\eta)}}{\Gamma(\alpha)}\int_a^x \frac{\tau^
{\rho\eta+\rho-1}\,f(\tau)}{(x^\rho-\tau^\rho)^{1-\alpha}}  d\tau  \label{eq:ek1}
\end{align}
 
for $x > a\geq 0, Re(\alpha)> 0$ and 

\begin{align}
(D^\alpha_{a+; \,\rho, \,\eta} f)(x) = x^{-\rho\eta}\bigg(\frac{1}{\rho\,x^{\rho-1}}\,\frac{d}{dx}\bigg)^n x^{\rho(n+\eta)}\, \Big(I^{n-\alpha}_{a+; \,\rho, \,\eta+\alpha} f\Big)(x) 
\label{eq:ek2}
\end{align}
for $x > a, Re(\alpha)\geq 0, \rho >0$, respectively. When $\rho =2, \, a=0$, the operators are called \emph{Erd\'{e}lyi-Kober} operators. When $\rho =1, \, a=0$, they are called \emph{Kober-Erd\'{e}lyi} or \emph{Kober operators} \cite[p.105]{key-8}. Extensive treatment of these operators can be found especially in \cite{key-Y}.
Further generalizations can be found, for example in \cite{key-8,key-9}. 

In \cite{key-0}, the author introduces an another generalization to the Riemann-Liouville and the Hadamard fractional integral and also provided existence results and semi-group properties. In \cite{udita2}, the author showed that both the Riemann-Liouville and the Hadamard fractional derivatives can be represented by a single fractional derivative operator. These new operators have been defined on the following extended Lebesgue space.

\section{Generalized Fractional Integration and Differentiation}

As in \cite{key-3}, consider the space $\textit{X}^p_c(a,b) \; (c\in \mathbb{R}, \, 1 \leq p \leq \infty)$ of those complex-valued Lebesgue measurable functions $f$ on $[a, b]$ for which $\|f\|_{\textit{X}^p_c} < \infty$, where the norm is defined by
\begin{equation}\label{eq:dff1}
\|f\|_{\textit{X}^p_c} =\Big(\int^b_a |t^c f(t)|^p \frac{dt}{t}\Big)^{1/p} < \infty, \quad (1 \leq p < \infty,\, c \in \mathbb{R}) 
\end{equation}
\noindent and for the case $p=\infty$,
\begin{equation} \label{eq:dff2}
\|f\|_{\textit{X}^\infty_c} = \text{ess sup}_{a \leq t \leq b} [t^c|f(t)|],  \quad ( c \in \mathbb{R}).
\end{equation}
In particular, when $c = 1/p \; (1 \leq p \leq \infty),$ the space $\textit{X}^p_c$ coincides with the classical $\textit{L}^p(a,b)$-space with
\begin{align} \label{eq:dff3}
&\|f\|_p =\Big(\int^b_a |f(t)|^p \, dt \Big)^{1/p} < \infty \quad (1 \leq p < \infty), \\
&\|f\|_\infty = \text{ess sup}_{a \leq t \leq b} |f(t)|   \quad ( c \in \mathbb{R}).
\end{align} 

Here we give the generalized forms of fractional integrals introduced in \cite{key-0} with a slight modification in the notations. 
\begin{definition}(Generalized Fractional Integrals)\\ 
Let $\Omega = [a,b] \, (-\infty <a < b < \infty)$ be a finite interval on the real axis $\mathbb{R}$. The generalized \emph{left-sided} fractional integral ${}^\rho I^\alpha_{a+}f$ of order $\alpha \in \mathbb{C} \; (Re(\alpha) > 0)$ of $f \in \textit{X}^p_c(a,b)$ is defined by
\begin{equation}
\big({}^\rho \mathcal{I}^\alpha_{a+}f\big)(x) = \frac{\rho^{1- \alpha }}{\Gamma({\alpha})} \int^x_a \frac{\tau^{\rho-1} f(\tau) }{(x^\rho - \tau^\rho)^{1-\alpha}}\, d\tau 
\label{eq:df1}
\end{equation}
for $x > a$ and $Re(\alpha) > 0$. The \emph{right-sided} generalized fractional integral
${}^\rho I^\alpha_{b-}f$ is defined by
\begin{equation}
\big({}^\rho \mathcal{I}^\alpha_{b-}f\big)(x) = \frac{\rho^{1- \alpha }}{\Gamma({\alpha})} \int^b_x \frac{\tau^{\rho-1} f(\tau) }{(\tau^\rho - x^\rho)^{1-\alpha}}\, d\tau
\label{eq:df2}
\end{equation}  
for $x < b$ and $Re(\alpha) > 0$. When $b=\infty$, the generalized fractional integral is called a Liouville-type integral.
\end{definition}

\begin{remark}
It can be seen that these operators are not equivalent to Erd\'{e}lyi-Kober operators due to the presence of a factor $\rho^{-\alpha}$, which is essential in the case of Hadamard operators. 
\end{remark}

Now consider the generalized fractional derivatives introduced in \cite{udita2};

\begin{definition}(Generalized Fractional Derivatives)\\
Let $\alpha \in \mathbb{C}, Re(\alpha) \geq 0, n=\ipa$ and $\rho >0.$ The generalized fractional derivatives, corresponding to the generalized fractional integrals (\ref{eq:df1}) and (\ref{eq:df2}), are defined, for $0 \leq a < x < b \leq \infty$, by
\begin{align}
\big({}^\rho \mathcal{D}^\alpha_{a+}f\big)(x) &= \bigg(x^{1-\rho} \,\frac{d}{dx}\bigg)^n\,\, \big({}^\rho \mathcal{I}^{n-\alpha}_{a+}f\big)(x)\nonumber\\
 &= \frac{\rho^{\alpha-n+1 }}{\Gamma({n-\alpha})} \, \bigg(x^{1-\rho} \,\frac{d}{dx}\bigg)^n \int^x_a \frac{\tau^{\rho-1} f(\tau) }{(x^\rho - \tau^\rho)^{\alpha-n+1}}\, d\tau
\label{eq:gd1}
\end{align}
and
\begin{align}
\big({}^\rho \mathcal{D}^\alpha_{b-}f\big)(x) &= \bigg(-x^{1-\rho} \,\frac{d}{dx}\bigg)^n\,\, \big({}^\rho \mathcal{I}^{n-\alpha}_{b-}f\big)(x)\nonumber\\
 &= \frac{\rho^{\alpha-n+1 }}{\Gamma({n-\alpha})} \, \bigg(-x^{1-\rho} \,\frac{d}{dx}\bigg)^n \int^b_x \frac{\tau^{\rho-1} f(\tau) }{(\tau^\rho - x^\rho)^{\alpha-n+1}}\, d\tau,
\label{eq:gd2}
\end{align}
respectively, if the integrals exist. When $b=\infty$, the generalized fractional derivative is called a Weyl or Liouville-type derivative \cite{udita2}.
\end{definition}

Next we give several important theorems related to these operators, without proofs. Interested reader may find them in \cite{udita2}. 


The following theorem gives the relations of the generalized fractional derivatives to that of Riemann-Liouville and Hadamard. For simplicity, we give only the \emph{left-sided} versions here. The similar results exist for the \emph{right-sided} operators as well.

\begin{theorem}\emph{\cite{udita2}}
Let $\alpha \in \mathbb{C},\, Re(\alpha) \geq 0,\, n=\lceil Re(\alpha) \rceil$ and $\rho >0.$ Then, for $x > a$,
\begin{align}
\noindent	&1. \lim_{\rho \rightarrow 1}\big({}^\rho \mathcal{I}^\alpha_{a+}f\big)(x) =\frac{1}{\Gamma(\alpha)}\int_a^x (x - \tau)^{\alpha -1} f(\tau) d\tau, \label{eq:th5} & & & &\\
	&2. \lim_{\rho \rightarrow 0^+}\big({}^\rho \mathcal{I}^\alpha_{a+}f\big)(x)=\frac{1}{\Gamma(\alpha)}\int_a^x \Bigg(\log\frac{x}{\tau}\Bigg)^{\alpha -1} f(\tau)\frac{d\tau}{\tau}, \label{eq:th2} & & & &\\
	&3. \lim_{\rho \rightarrow 1}\big({}^\rho \mathcal{D}^\alpha_{a+}f\big)(x)=\bigg(\frac{d}{dx}\bigg)^n \,\frac{1}{\Gamma({n-\alpha})} \, \int^x_a \frac{ f(\tau) }{(x - \tau)^{\alpha-n+1}}\, d\tau,\label{eq:th3} & & & &\\
	&4. \lim_{\rho \rightarrow 0^+}\big({}^\rho \mathcal{D}^\alpha_{a+}f\big)(x)=\frac{1}{\Gamma(n-\alpha)}\Bigg(x\frac{d}{dx}\Bigg)^n \int_a^x \Bigg(\log\frac{x}{\tau}\Bigg)^{n-\alpha +1} f(\tau)\frac{d\tau}{\tau}.\label{eq:th4}& & & &
\end{align}
\end{theorem}
\begin{proof}
We reproduce the proof here for completeness. The equations (\ref{eq:th5}) and (\ref{eq:th3}) follow from taking limits when $\rho \rightarrow 1$, while (\ref{eq:th2}) follows from the L'H\^ospital rule by noticing that \cite{key-0}
\begin{align}
\lim_{\rho \rightarrow 0^+}\frac{\rho^{1- \alpha }}{\Gamma({\alpha})} \int^x_a \frac{ f(\tau)\tau^{\rho-1} }{(x^\rho - \tau^\rho)^{1-\alpha}}\, d\tau 
  &= \frac{1}{\Gamma({\alpha})}\int^x_a \lim_{\rho \rightarrow 0^+} f(\tau)\tau^{\rho-1} \Bigg(\frac{x^\rho - \tau^\rho}{\rho}\Bigg)^{\alpha -1}d\tau \nonumber\\
	&= \frac{1}{\Gamma({\alpha})}\int^x_a  f(\tau)\lim_{\rho \rightarrow 0^+}\Bigg(\frac{x^\rho - \tau^\rho}{\rho}\Bigg)^{\alpha -1}\frac{d\tau}{\tau} \nonumber\\
	&= \frac{1}{\Gamma({\alpha})}\int^x_a  f(\tau)\Bigg(\lim_{\rho \rightarrow 0^+}\frac{x^\rho - \tau^\rho}{\rho}\Bigg)^{\alpha -1}\frac{d\tau}{\tau} \label{eql1}\\
	&= \frac{1}{\Gamma(\alpha)}\int_a^x \Bigg(\log\frac{x}{\tau}\Bigg)^{\alpha -1} f(\tau)\frac{d\tau}{\tau}, \nonumber
\end{align}
taking into account the principal branch of the logarithm in equation (\ref{eql1}). The proof of (\ref{eq:th4}) is similar.
\end{proof}

\begin{remark}
Note that the equations (\ref{eq:th5}) and (\ref{eq:th3}) are related to the Riemann-Liouville operators, while the equations (\ref{eq:th2}) and (\ref{eq:th4}) are related to the Hadamard operators. The results similar to (\ref{eq:th5}) and (\ref{eq:th3}) can also be derived from Erd\'{e}lyi-Kober operators, though it is not possible to obtain equivalent results for (\ref{eq:th2}) and (\ref{eq:th4}) due to the absence of the factor $\rho^{\alpha}$, which is necessary in the case of the Hadamard-type operators.
\end{remark}

Next is the inverse property. 

\begin{theorem}\label{th:inv}\emph{\cite{udita2}}
Let $0 <\alpha <1$, and $f \in \textit{X}^p_c(a,b)$. Then, for $a>0,\, \rho >0$,
\begin{equation}
  \Big({}^\rho \mathcal{D}^\alpha_{a+}\,{}^\rho \mathcal{I}^\alpha_{a+}\Big)f(x) = f(x).
\end{equation}
\end{theorem}

Compositions between the operators of generalized fractional differentiation and generalized fractional integration are given by the following theorem.
\begin{theorem}\emph{\cite{udita2}}
Let $\alpha,\, \beta \in \mathbb{C}$ be such that $0 < Re(\alpha) < Re(\beta) < 1.$ If $0 < a < b < \infty$ and $ 1 \leq p \leq \infty$, then, for $f \in L^p(a,b), \,\, \rho >0$,
\begin{equation}
  {}^\rho \mathcal{D}^\alpha_{a+}\,{}^\rho \mathcal{I}^\beta_{a+}\,f = {}^\rho \mathcal{I}^{\beta - \alpha}_{a+}\,f \;\;\text{and}\;\;{}^\rho \mathcal{D}^\alpha_{b-}\,{}^\rho \mathcal{I}^\beta_{b-}\,f = {}^\rho \mathcal{I}^{\beta - \alpha}_{b-}\,f. \nonumber
\end{equation}
\end{theorem}

In \cite{udita2}, the author derived the formula for the generalized derivative (\ref{eq:gd1}) of the power function $f(x) = x^\nu$, where $\nu \in \mathbb{R}$,  for $\alpha \in \mathbb{R}^+, \, 0 < \alpha < 1$ and $a=0$, given by, 
\begin{align}
    {}^\rho D^\alpha_{0+} x^\nu &= \frac{\Gamma\Big(1+\frac{\nu}{\rho}\Big)\,\rho^{\alpha-1}}{\Gamma\Big(1+\frac{\nu}{\rho}-\alpha\Big)}\, x^{\nu -\alpha\rho} 
\label{eq:gen}    
\end{align}
for $\rho > 0$. When $\rho = 1$, this agrees with the standard results obtained for the Riemann-Liouville fractional derivative (\ref{eq:lRLD}) \cite{key-8,key-2,key-9}.
Interestingly enough, for $\alpha = 1,\, \rho = 1$, we obtain ${}^1D^1_{0+} x^\nu = \nu\,x^{\nu -1}$, as one would expect.

\begin{figure}[h]
	\centering
	  \subfigure[$\nu$= 2.0]{\includegraphics[width=2.5in, height=1.8in, angle=0, clip]{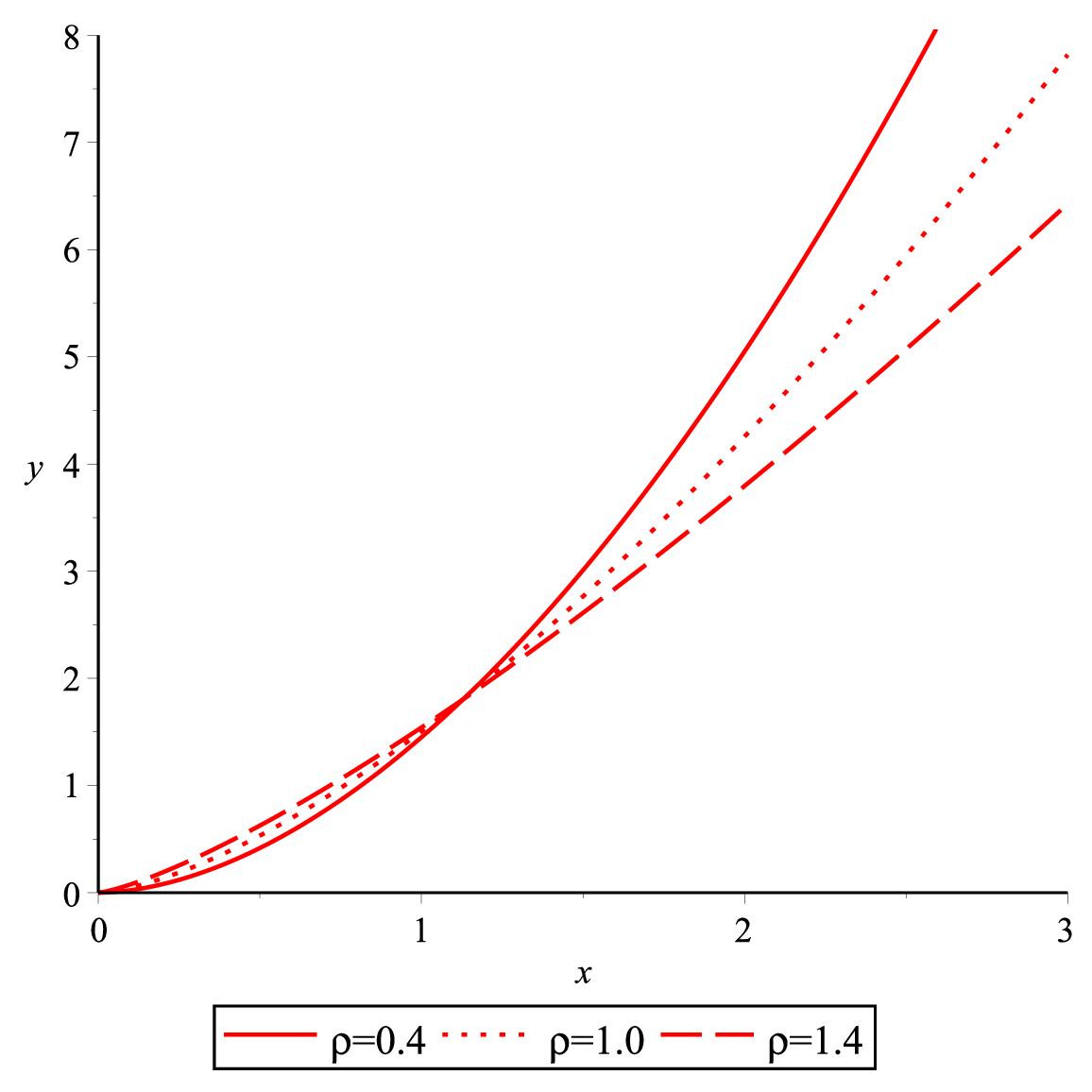}}
	  \hspace{0.01in}
	  \subfigure[$\nu$= 0.5]{\includegraphics[width=2.6in, height=1.8in, angle=0, clip]{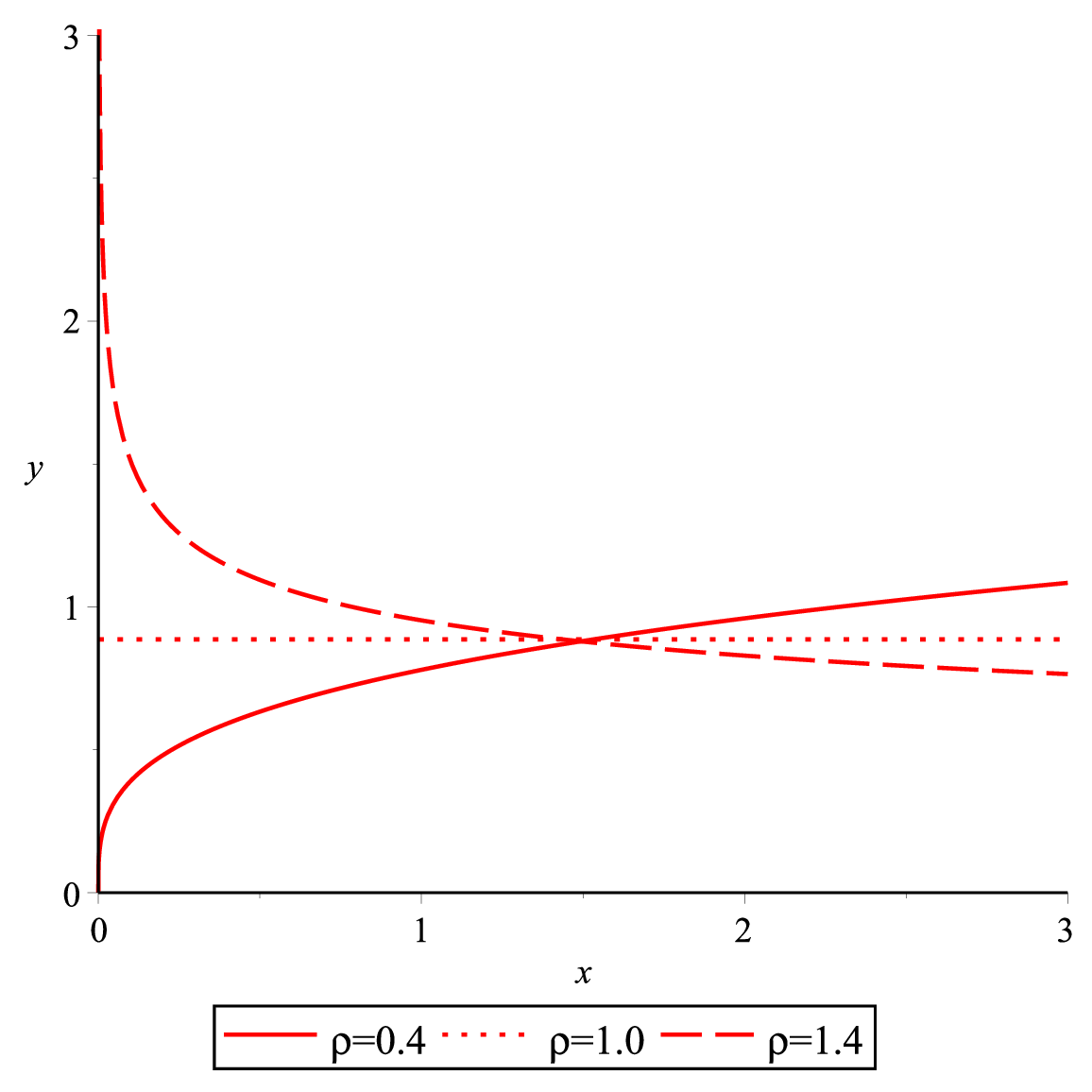}}\\
	  
	  \subfigure[$\nu$= 2.0]{\includegraphics[width=2.5in, height=1.8in, angle=0, clip]{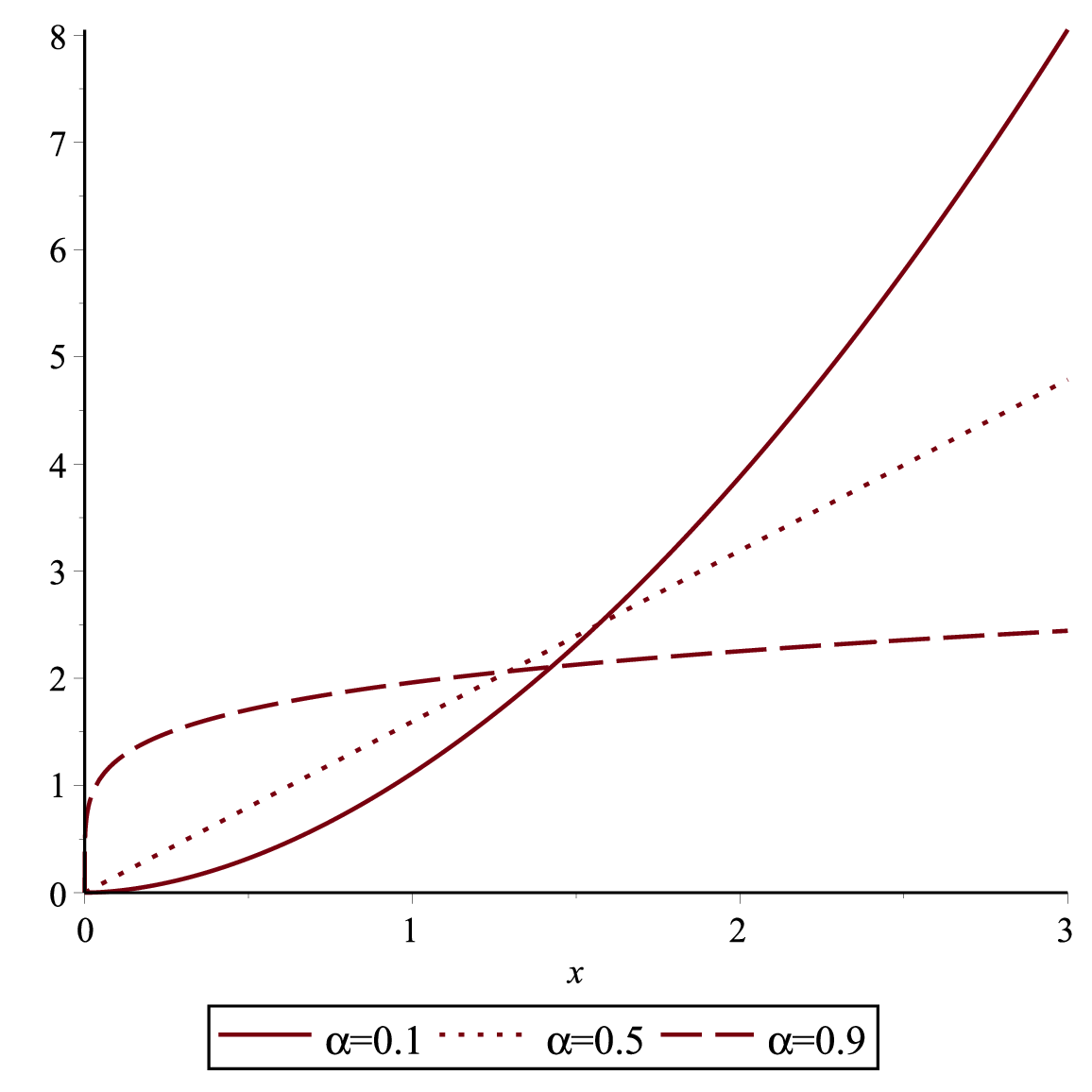}}
	  \hspace{0.01in}
	  \subfigure[$\alpha$= 0.5]{\includegraphics[width=2.5in, height=1.8in, angle=0, clip]{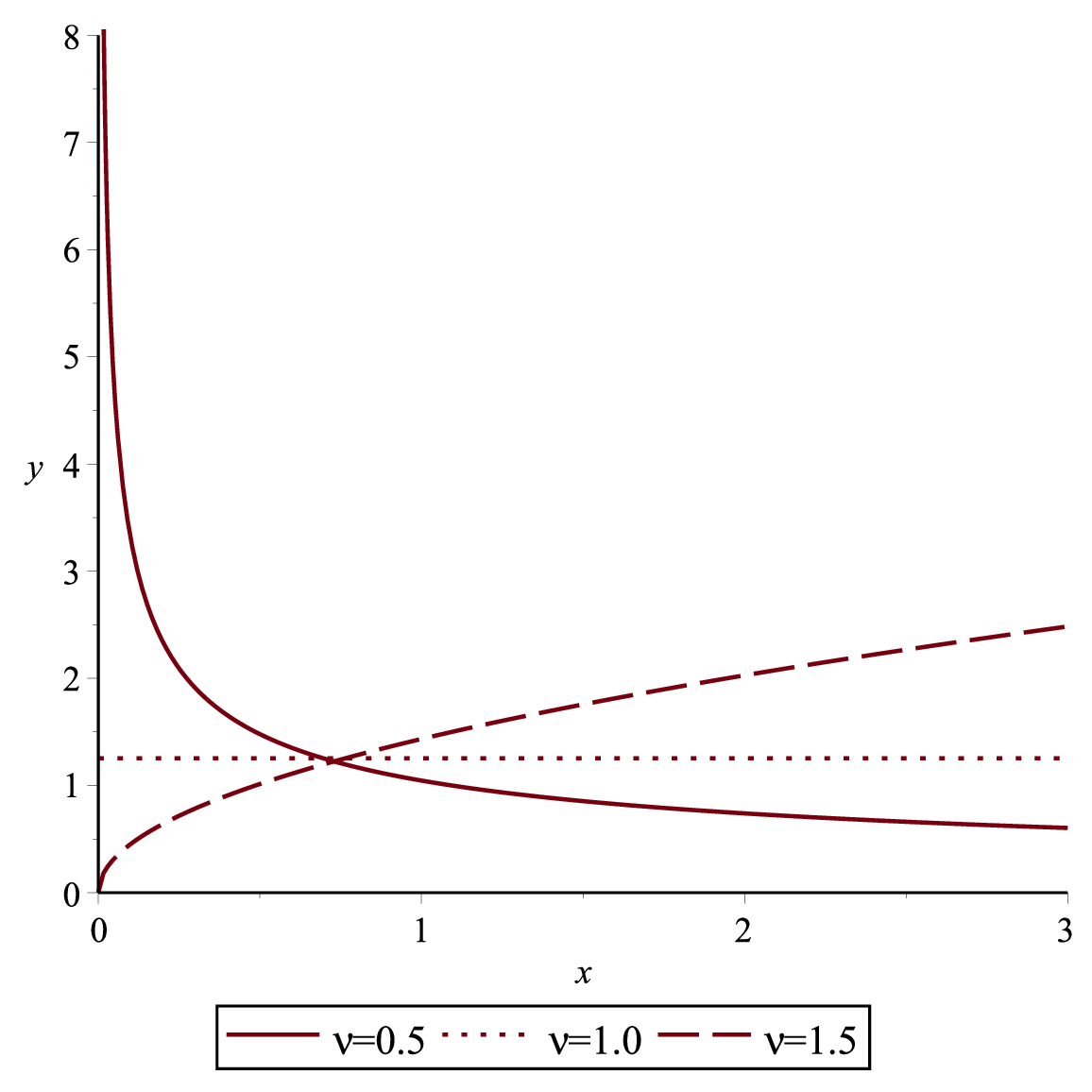}}	  
	\caption{Generalized fractional derivative of the power function $x^\nu$ for $\rho = 0.4,\, 1.0,\, 1.4$, $\nu = 0.5, \, 2.0$ and $\alpha = 0.1, \, 0.5, \, 0.9$ \cite{udita2}}        
	\label{fig:FD-1}
\end{figure}

Figure \ref{fig:FD-1} shows the comparison results of the generalized derivative (\ref{eq:gen}) for different values of $\rho \, (\mbox{for} \, \alpha=1/2)$, $\nu \, (\mbox{for} \, \alpha=1/2)$ and $\alpha$. 


In the next section we give the main results of this paper, the Mellin transforms 
of the generalized fractional operators. 
\section{Mellin Transforms of Generalized Fractional Operators}
According to Flajolet et al. \cite{key-15}, Hjalmar Mellin\footnote{H. Mellin was a Finnish and his adviser was Gosta Mittag-Leffler, a Swedish mathematician. Later he also worked with. Kurt Weierstrass \cite{key-22}.}(1854-1933) gave his name to the \emph{Mellin transform} that associates to a function $f(x)$ defined over the positive reals, the complex function $\mathcal{M}[f]$ \cite{key-16}. It is closely related to the \emph{Laplace} and \emph{Fourier} transforms. 

We start by recalling the important properties of the Mellin transform. The domain of definition is an open strip, $<a,b>$, say, of complex numbers $s= \sigma + it$ such that $0 \leq a < \sigma < b$. Here we adopt the definitions and properties mentioned in \cite{key-15} with some minor modifications to the notations. 
\begin{definition}(Mellin transform) 
Let $f(x)$ be locally Lebesgue integrable over $(0, \infty)$. The \emph{Mellin transform} of $f(x)$ is defined by
$$
 \mathcal{M}[f](s) = \int_0^\infty x^{s-1}f(x) \,dx. \label{eq:MT}
$$
The largest open strip $<a, b>$ in which the integral converges is called the \emph{fundamental strip}.
\end{definition}
Following theorem will be of great importance in applications. 
\begin{theorem}[(\cite{key-15}, Theorem 1)] Let $f(x)$ be a function whose transform admits the fundamental strip $<a, b>$. Let $\rho$ be a nonzero real number, and $\mu$, $\nu$ be positive reals. Then,
\begin{align*}
 &1. \;\mathcal{M}\bigg[\sum_k \lambda_k f(\mu_k x)\bigg](s) = \Bigg( \sum_{k\in I} \frac{\lambda_k}{\mu_k^s} \Bigg) \mathcal{M}[f](s), \quad I \;\text{finite}, \,\lambda_k >0, \,s \in <a, b>& & & \\
 &2. \;\mathcal{M}\big[x^\nu f(x)\big](s) = \mathcal{M}[f](s + \nu) \quad  s \in <a, b>& & & \\
 &3. \;\mathcal{M}\big[f(x^\rho)\big](s) = \frac{1}{\rho} \mathcal{M}[f]\bigg(\frac{s}{\rho}\bigg), \quad s \in <\rho a, \rho b>& & & \\
 &4. \;\mathcal{M}\Big[\frac{d}{dx}f(x)\Big](s) = (1-s)\, \mathcal{M}[f](s - 1) & & & \\
 &5. \;\mathcal{M}\Big[x\,\frac{d}{dx}f(x)\Big](s) = -s\, \mathcal{M}[f](s) & & & \\
 &6. \;\mathcal{M}\Big[\int_0^x f(t) dt \Big](s) = - \frac{1}{s} \, \mathcal{M}[f](s + 1) & & & 
\end{align*}
\end{theorem}
Next is the inversion theorem for \emph{Mellin transforms}.  
\begin{theorem}[(\cite{key-15}, Theorem 2)] 
Let $f(x)$ be integrable with fundamental strip $<a, b>$. If $c$ is such that $a < c < b$, and $\mathcal{M}[f](c+it)$ integrable, then the equality,
$$
\frac{1}{2\pi i}\, \int_{c \,- \,i\infty}^{c\, + \,i\infty} \mathcal{M}[f](s)\, x^{-s}\, ds = f(x)
$$
holds almost everywhere. Moreover, if $f(x)$ is continuous, then equality holds everywhere on $(0, \infty)$. 
\end{theorem}
Following is the first important result of the paper, i.e., the Mellin transforms of the generalized fractional integrals. We give both the \emph{left-sided} and \emph{right-sided} versions of the results here.
\begin{lemma}
Let $\alpha \in \mathbb{C},\, Re(\alpha) > 0,$ and $\rho >0.$ Then, 
\begin{align}
  &\mathcal{M}\bigg({}^\rho \mathcal{I}^\alpha_{a+}f\bigg)(s) = \frac{\Gamma\big(1-\frac{s}{\rho}-\alpha\big)}{\Gamma\big(1-\frac{s}{\rho}\big)\,\rho^\alpha}\, \mathcal{M}f(s + \alpha\rho), \quad Re(s/\rho + \alpha) < 1, \, x > a, & & & \label{eq:MT1}\\
  &\mathcal{M}\bigg({}^\rho \mathcal{I}^\alpha_{b-}f\bigg)(s) = \frac{\Gamma\big(\frac{s}{\rho}\big)}{\Gamma\big(\frac{s}{\rho} + \alpha\big)\,\rho^\alpha}\, \mathcal{M}f(s + \alpha\rho), \quad Re(s/\rho) > 0, \, x < b, & & & \label{eq:MT2} 
\end{align}
for $f \in \textit{X}^1_{s + \alpha\rho}(\mathbb{R^+})$, if $\mathcal{M}f(s + \alpha\rho)$ exists for $s\in \mathbb{C}$. 
\end{lemma}
\begin{proof}
We use the Fubini's theorem and the Dirichlet technique here. By the Definition \ref{eq:MT} and Equation (\ref{eq:df1}), we have
\begin{align*}
\mathcal{M}\bigg({}^\rho \mathcal{I}^\alpha_{a+}f\bigg)(s) &= \int_0^\infty x^{s-1} \, \frac{\rho^{1-\alpha}}{\Gamma(\alpha)}\int_a^x \, (x^\rho - \tau^\rho)^{\alpha -1} \tau^{\rho - 1} f(\tau) \,d\tau \,dx, \\
 &= \frac{\rho^{1-\alpha}}{\Gamma(\alpha)}\int_0^\infty \tau^{\rho - 1} f(\tau) \int_\tau^\infty x^{s-1} (x^\rho - \tau^\rho)^{\alpha -1} dx\,d\tau\\
 &= \frac{\rho^{-\alpha}}{\Gamma(\alpha)}\int_0^\infty \tau^{s+\alpha\rho -1} f(\tau) \int_0^1 u^{-\frac{s}{\rho} -\alpha} (1-u)^{\alpha - 1} du\,d\tau \\
 & =\frac{\rho^{-\alpha}\Gamma\big(1-\frac{s}{\rho} -\alpha\big)}{\Gamma\big(1-\frac{s}{\rho}\big)} \int_0^\infty \tau^{s+\alpha\rho -1} f(\tau) \\ &=\frac{\Gamma\big(1-\frac{s}{\rho} -\alpha\big)}{\Gamma\big(1-\frac{s}{\rho}\big) \rho^{\alpha}}\, \mathcal{M}[f] (s + \alpha\rho) \quad \text{for} \,\,Re(s/\rho + \alpha) < 1.
\end{align*}
after using the change of variable $u=(\tau/x)^\rho$, and properties of the Beta function. This proves (\ref{eq:MT1}). The proof of (\ref{eq:MT2}) is similar. 
\end{proof}
\begin{remark}
The two transforms above confirm Lemma 2.15 of \cite{key-8} for $\rho = 1$. For the case, when $\rho \rightarrow 0^+$, consider the quotient expansion of two $Gamma$ functions at infinity given by \cite{key-4},
\begin{equation}
   \frac{\Gamma(z+a)}{\Gamma(z+b)} = z^{a-b}\, \Big[1+O\bigg(\frac{1}{z}\bigg)\Big],   \quad |arg(z+a)| < \pi; \,\, |z| \rightarrow \infty. \label{eq:qGamma}
\end{equation} 
Then, it is clear that, 
\begin{align*}
\lim_{\rho \rightarrow 0^+}\big({}^\rho \mathcal{I}^\alpha_{a+}f\big)(x) &= \lim_{\rho \rightarrow 0^+} \frac{\Gamma\big(1-\frac{s}{\rho} -\alpha\big)}{\Gamma\big(1-\frac{s}{\rho}\big) \rho^{\alpha}}\, \mathcal{M}[f] (s + \alpha\rho)\\
&= (-s)^{-\alpha}\,\mathcal{M}[f](s), \quad Re(s) < 0,
\end{align*}
where $(-s)^{-\alpha}=e^{-\alpha\log(-s)}$, taking into account the principal branch of the logarithm. This confirms Lemma 2.38(a) of \cite{key-8}.  
\end{remark}
To prove the next result we use the Mellin transform of $m^{th}$ derivative of a m-times differentiable function given by, 
\begin{lemma}[(\cite{key-8}, p.21)]
Let $\varphi \in \mathbf{C}^m(\mathbb{R}^+)$, $\mathcal{M}[\varphi(t)](s-m)$ and $\mathcal{M} [D^m\varphi(t)](s)$ exist, and $lim_{t \rightarrow 0^+} \big[t^{s-k-1}\varphi^{(m-k-1)}(t)\big]$ and $lim_{t \rightarrow +\infty} \big[t^{s-k-1}\varphi^{(m-k-1)}(t)\big]$ are finite for $k=0, 1, \cdots , m-1$, $m \in \mathbb{N}$, then
\begin{align}
 \mathcal{M} \big[D^m\varphi(t)\big](s) = &\frac{\Gamma(1+m-s)}{\Gamma(1-s)}\, \mathcal{M}[\varphi](s-m)\nonumber\\
 &\hspace{.5cm}+\sum_{k=0}^{m-1}\frac{\Gamma(1+k-s)}{\Gamma(1-s)}\,\big[x^{s-k-1} \varphi^{(m-k-1)}(x)\big]^\infty_0 \label{eq:MD1}
\end{align}
\begin{align}
 \mathcal{M} \big[D^m\varphi(t)\big](s)= (-1)^m &\frac{\Gamma(s)}{\Gamma(s-m)}\, \mathcal{M}[\varphi](s-m)\nonumber\\
 &+\sum_{k=0}^{m-1} (-1)^k \frac{\Gamma(s)}{\Gamma(s-k)}\,\big[x^{s-k-1} \varphi^{(m-k-1)}(x)\big]^\infty_0 \label{eq:MD2}
\end{align}
\end{lemma}
The next result is the Mellin transforms of the generalized fractional derivatives. For simplicity we consider only the case $0 <\alpha < 1$ here. In this case, $n = \lceil \alpha \rceil = 1$. 
\begin{theorem}
Let $\alpha \in \mathbb{C},\, Re(\alpha) > 0 ,\, s \in \mathbb{C},\, \rho >0$ and $f(x) \in \textit{X}^1_{s-\alpha\rho}(\mathbb{R}^+)$. Also assume $f(x)$ satisfies the following conditions:
\[
\lim_{x \rightarrow 0^+} \,x^{s-\rho}\Big({}^\rho \mathcal{I}^{1-\alpha}_{a+}f\Big)(x) = 
\lim_{x \rightarrow \infty} \,x^{s-\rho}\Big({}^\rho \mathcal{I}^{1-\alpha}_{a+}f\Big)(x) = 0,
\]
and 
\[
\lim_{x \rightarrow 0^+} \,x^{s-\rho}\Big({}^\rho \mathcal{I}^{1-\alpha}_{b-}f\Big)(x) = 
\lim_{x \rightarrow \infty} \,x^{s-\rho}\Big({}^\rho \mathcal{I}^{1-\alpha}_{b-}f\Big)(x) = 0,
\]
respectively. Then, 
\begin{align}
  &\mathcal{M}\bigg({}^\rho \mathcal{D}^\alpha_{a+}f\bigg)(s) =  \frac{\rho^\alpha\,\Gamma\big(1-\frac{s}{\rho}+\alpha\big)}{\Gamma\big(1-\frac{s}{\rho}\big)}\, \mathcal{M}[f](s - \alpha\rho), \quad Re(s/\rho) < 1, \, x > a \geq 0, & &  \label{eq:MTD1}\\
  &\mathcal{M}\bigg({}^\rho \mathcal{D}^\alpha_{b-}f\bigg)(s) =  \frac{\rho^\alpha\,\Gamma\big(\frac{s}{\rho}\big)}{\Gamma\big(\frac{s}{\rho} - \alpha\big)}\, \mathcal{M}[f](s - \alpha\rho), \quad Re(s/\rho - \alpha) > 0, \, x < b \leq \infty, & &  \label{eq:MTD2} 
\end{align}
respectively. 
\end{theorem}
\begin{proof}
By Definition \ref{eq:MT} and Equation (\ref{eq:gd1}), we have
\begin{align}
\mathcal{M}\bigg({}^\rho \mathcal{D}^\alpha_{a+}f\bigg)(s) &= \frac{\rho^{\alpha}}{\Gamma(1-\alpha)}\int_0^\infty x^{s-1} \, \Big(x^{1-\rho}\frac{d}{dx}\Big) \int_a^x \frac{\tau^{\rho - 1}}{(x^\rho - \tau^\rho)^{\alpha } }  f(\tau) \,d\tau \,dx, \label{eq:prod1}\\
 &=\int_0^\infty x^{s-\rho} \, \frac{d}{dx} \bigg({}^\rho \mathcal{I}^{1-\alpha}_{a+}f\bigg)(x) dx \nonumber\\ &=\frac{\rho^\alpha\,\Gamma\big(1-\frac{s}{\rho}+\alpha\big)}{\Gamma\big(1-\frac{s}{\rho}\big)} \, \mathcal{M}[f] (s - \alpha\rho) + x^{s-\rho}\Big({}^\rho \mathcal{I}^{1-\alpha}_{a+}f\Big)(x)\Big|^\infty_a \nonumber
\end{align}
for $Re(s/\rho) < 1$, using (\ref{eq:MD1}) with $m=1$ and (\ref{eq:df1}). This establishes (\ref{eq:MTD1}). The proof of (\ref{eq:MTD2}) is similar. 
\end{proof}
\begin{remark}
The Mellin transforms of Riemann-Liouville and Hadamard fractional derivatives follow easily, i.e., the transforms above confirm Lemma 2.16 of \cite{key-8} for $\rho = 1, \, b=\infty$ and Lemma 2.39 of \cite{key-8} for the case when $\rho \rightarrow 0^+$ in view of (\ref{eq:qGamma}).  
\end{remark}

\begin{remark}
In many cases we only considered \emph{left-sided} derivatives and integrals. But the \emph{right-sided} operators can be treated similarly. 
\end{remark}

\begin{remark}
For $ Re(\alpha) > 1$, we need to replace $x^{1-\rho} \frac{d}{dx}$ in (\ref{eq:prod1}) by $(x^{1-\rho} \frac{d}{dx})^{\lceil Re(\alpha) \rceil}$. In this case the equation becomes very complicated and introduces several interesting combinatorial problems. There is a well-developed branch of mathematics called \emph{Umbral Calculus} \cite{key-24,key-29}, which consider the case when $1-\rho \in \mathbb{N}$. Bucchianico and Lobe \cite{key-25} provide an extensive survey of the subject. We discuss some of them here referring interested reader to the works of Rota, Robinson and Roman \cite{key-30,key-26,key-27,key-28,key-29}. First consider the $ generalized  \,\, \delta$-derivative defined by,
\end{remark}

\begin{definition}[]
Let $r \in \mathbb{R}$ and $m \in \mathbb{N}_0 = \{0, 1, 2, 3, \ldots \}$. The generalized $\delta_{r,m}$-differential operator is defined by
\begin{equation}
\delta_{r, m} := x^r \, \frac{d^m}{dx^m}. \label{eq:dd}
\end{equation}
\end{definition}
The $n^{th}$ degree $\delta_{r, m}$ - derivative operator has interesting properties. When $r=m=1$ it generates the \emph{Stirling numbers of the second kind, $S(n, k)$} \cite{key-23}, \,Sloane's A008277 \cite{key-21}, while $r=2$ and $m=1$ generate the \emph{unsigned Lah numbers, $L(n, k)$},\, Sloane's A008297, A105278 \cite{key-21}. They have been well-studied and appear frequently in literature of combinatorics \cite{key-20}. These can also be looked at as the number of ways to place $j$ non-attacking rooks on a Ferrer's-board with certain properties based on $r$ \cite{key-19}.\\ 

Now, let $c_{r, m; \,j}^n, \; j = 1,\, 2,\, 3, \,\cdots, \,n$; $n \in \mathbb{N}, \; r \in \mathbb{R}$  be the $j^{th}$ coefficient of the expansion of $\delta_{r, m}^n$, in a basis, $\mathcal{B} = \Big\{x^{n(k-1)+1}\frac{d^m}{dx^m},\, x^{n(k-1)+2}\frac{d^{m+1}}{dx^{m+1}},\, \cdots,\, x^{nk}\frac{d^{n+m}}{dx^{n+m}}\Big\}$, i.e., 
\begin{equation}
\Bigg(x^k \frac{d^m}{dx^m}\Bigg)^n = c_{r,m; \,1}^n \,x^{n(k-1)+1}\frac{d^m}{dx^m} + c_{r, m; \,2}^n\,x^{n(k-1)+2}\frac{d^{m+1}}{dx^{m+1}} + \cdots + c_{r, m; \,n}^n \,x^{nk}\frac{d^{n+m}}{dx^{n+m}}
\label{eq:d-exp}
\end{equation}
Table 1 lists the cases for $r=1$ and $r=2$ when $m=1$ in triangular settings,\\ 

\hspace*{5.85cm} 1 \hspace{3.95cm} 1 \\ 
\hspace*{5.95cm} 1 \hspace{.5cm} 1 \hspace{3.25cm}2 \hspace{.45cm}1 \\
\hspace*{5.6cm} 1\hspace{.5cm} 3 \hspace{.5cm}1\hspace{2.7cm}6\hspace{.5cm}6\hspace{.5cm}1\\
\hspace*{5.4cm}1\hspace{.5cm}7\hspace{.6cm}6\hspace{.5cm}1\hspace{1.95cm}24\hspace{.3cm}36\hspace{.4cm}12\hspace{.4cm}1\\
\hspace*{5.1cm}1\hspace{.35cm}15\hspace{.4cm}25\hspace{.4cm}10\hspace{.4cm}1\hspace{1.2cm}120\hspace{.1cm}\hspace{.1cm}240\hspace{.1cm}120\hspace{.25cm}20\hspace{.45cm}1\\
\hspace*{4.8cm}.\hspace{.6cm}.\hspace{.6cm}.\hspace{.6cm}.\hspace{.6cm}.\hspace{.6cm}.\hspace{.8cm}.\hspace{.6cm}.\hspace{.6cm}.\hspace{.6cm}.\hspace{.6cm}.\hspace{.6cm}.\\\\
\hspace*{6.3cm}$r=1$\hspace{3.45cm}$r=2$\\\\
\hspace*{4.5cm}{Table 1: Coefficients of $n^{th}$ degree generalized $\delta$-derivative}\\

\normalsize 

The Stirling numbers of the second kind, $S(n, k)$ and unsigned Lah numbers, $L(n, k)$ are given by \cite{key-20}
\begin{equation}
  S(n, k) = \frac{1}{k!} \sum_{i=0}^k (-1)^i \, \binom{k}{i} \, \big(k-i\big)^n \quad ;0 \le k \le n \in \mathbb{N} \cup \{0\},
\end{equation}
and
\begin{equation}
  L(n, k) = \binom{n-1}{k-1} \frac{n!}{k!} \quad ;0 \le k \le n \in \mathbb{N} \cup \{0\},
\end{equation}  
respectively.

A closed-form formula exists for the case $r=2$, and asymptotic formula exists for $r=1$ \cite{key-17,key-18}. Johnson \cite{key-19} used analytical methods and Goldman \cite{key-17} used algebraic methods to obtain generating functions for $r=1, \, 2$ given by
\begin{equation}
  z(x, y) = \begin{cases}
    \exp\Big\{\frac{e^{xy}-1}{x}\Big\} & \quad \text{if $r = 1$},\\
    \exp\Big\{\frac{1}{x}\big[1-(r-1)xy\big]^{1/(1-r)} - \frac{1}{x}\Big\} & \quad \text{if $r \geq 2$}.\\
  \end{cases}
\end{equation}
It is interesting to note that the order of $x^r$ and $\frac{d^m}{dx^m}$ in the $\delta_{r, m}$ operator is of no importance as indicated in the following examples, as long as $r$ and $m$ are kept fixed. For example, we consider $\delta_{2, 1}$ and $\delta_{1, 2}$. The first five rows of the corresponding Pseudo polynomials are given below:
\begin{align*} 
&\dx{2} \\
2\dx{3}&+\ddx{4}{2}\\
&\hspace{-2cm}6\dx{4}+6\ddx{5}{2}+\ddx{6}{3}\\
24\dx{5}+36\ddx{6}{2}&+12\ddx{7}{3}+\ddx{8}{4}\\
&\hspace{-4.0cm}120\dx{6}+240\ddx{7}{2}+120\ddx{8}{3}+20\ddx{9}{4}+\ddx{10}{5}\\
\end{align*}
\hspace*{5.6cm}{Table 2: Pseudo polynomials of $\delta_{2,1}$                          

\begin{align*} 
&\xdx{2} \\
2\xdx{3}&+\ddx{2}{4}\\
&\hspace{-2cm}6\xdx{4}+6\ddx{2}{5}+\ddx{3}{6}\\
24\xdx{5}+36\ddx{2}{6}&+12\ddx{3}{7}+\ddx{4}{8}\\
&\hspace{-4.0cm}120\xdx{6}+240\ddx{2}{7}+120\ddx{3}{8}+20\ddx{4}{9}+\ddx{5}{10}
\end{align*}
\hspace*{5.2cm}{Table 3: Pseudo polynomials of $\delta_{1,2}$\\                         

These polynomials can be generated using the Maple\textsuperscript{\small{\textregistered}}codes given in the Algorithm \ref{alg} for $r=2$ and $m=1$. \normalsize

This interesting result holds for any $r \in\mathbb{N}=\{1, 2, 3, \cdots\}$ and will be given in the following result. In Theorem \ref{rec2}, we further show that the differential operators $x^r\frac{d^m}{dx^m}$ and $x^m\frac{d^r}{dx^r}$ produce the same sequence of coefficients as well.

\begin{theorem}\label{rec1}
Let $r\in\mathbb{N}=\{1, 2, 3, \cdots\}$. The generalized Stirling numbers generated by the operators $\delta_{r,1}$ and $\delta_{1, r}$ are coincide with each other and follow the recurrence relation given by
\begin{equation}\label{eq-rec1}
S(n,k)=\big[(n-1)(r-1)+k\big]S(n-1,k)+S(n-1,k-1), \;\;  0\leq k\leq n,
\end{equation}
with $S(0,0)=1$ and $S(n,0)=S(n,k)=0$ for $k\geq n+r\geq 0$. 
\end{theorem} 

It is clear from (\ref{eq-rec1}) that when $r=1$, the recurrence relation coincide with that of the \emph{Stirling numbers of the second kind} and when $r=2$, it coincide with that of the {unsigned Lah numbers}. 

\begin{proof}
We prove Theorem \ref{rec1} by induction on $n$ and $k$ for $\delta_{r,1}$. For $n=k=1$, we have $S(1,1)=S(0,1)+S(0,0)=1$. Now, for an inductive argument, assume the recurrence relation (\ref{eq-rec1}) holds for $n$ and $k$ and prove the relation holds for $n+1$ and $k+1$. First consider the case $S(n+1,1)$. The coefficient of $x^{r+(r-1)(n-1)}\frac{d}{dx}$ is $S(n,1)$ and applying the operator $x^r\frac{d}{dx}$, we see that $S(n+1,1)=(r+(r-1)(n-1))S(n,1)$. For $S(n+1,n+1)$, notice that $S(n,n+1)=0$ and $S(n,n)=1$ and by applying the differential operator $x^r\frac{d}{dx}$, we have $S(n+1,n+1)=S(n,n)$. For $k\notin\{1, n+1\}$, the coefficient $S(n+1,k)$ of $x^{r+(r-1)(n-1)+k}\frac{d^{k}}{dx^{k}}$ is obtained by applying $\delta_{r,1}$ to two therms in the $n^{th}$ step, namely, $S(n,k-1)x^{r+(r-1)(n-2)+k-1}\frac{d^{k-1}}{dx^{k-1}}$ and $S(n,k)x^{r+(r-1)(n-2)+k}\frac{d^k}{dx^k}$ and collecting proper terms. In the first product, we differentiate the second term, thus keep the coefficient fixed. From the second product, we obtain $S(n,k)(r+(r-1)(n-2)+k)$. This proves the equation (\ref{eq-rec1}) for $n$. 

Now, for the case $k$, we use a similar argument. It is clear that $S(k,k+1)=0$ and $S(k+1,k+1)=1$ for $k\in\mathbb{N}$. The coefficients $S(k+2,k+1), S(k+3,k+1), \ldots, S(n,k+1)$ can be obtained following the case of $n$. The case of $\delta_{1, r}$ is similar. This completes the proof.
\end{proof}

The Theorem \ref{rec1} can be extended to the cases of any $r, m \in \mathbb{N}$. Let us first consider, for example, the cases of $\delta_{2,3}$ and $\delta_{3, 2}$. In both the cases, it generates the same triangular array of sequences given in Table 4. This interesting result holds for any two positive integers and is given in Theorem \ref{rec2}.
\begin{table}[h]
\centering
\begin{tabular}{c|cccccccc}
\small{\backslashbox{n}{k}} &$c_{n,1}$ & $c_{n,2}$ & $c_{n,3}$ & $c_{n,4}$ & $c_{n,5}$ & $c_{n,6}$ & $c_{n,7}$ & $c_{n,8}$\\\hline
1 & 1 & &&&&&&\\
2 & 6 & 6 & 1\\
3 & 72 &168 & 96 & 18 & 1\\
4 & 1440 & 5760 & 6120 & 2520 & 456 & 36 & 1 \\
$\vdots$ & &&$\vdots$ & & &&&$\ddots$
\end{tabular}\\
\vspace{0.2in}
\begin{center}
Table 4. The sequence for $\delta_{2,3}$ and $\delta_{3,2}$ 
\end{center}
\end{table}


\subsection{Hidden Pascal Triangles}\label{pas}

Even though the terms in the sequence $\delta_{3,2}$ have no apparent resemblance to any familiar sequences, it is interesting to notice that there is a hidden Pascal Triangle buried inside those terms. The reason is that, in each differentiation step, it creates two terms, one with a degree one less than the degree of $x$ and one with a order one greater than that of $\frac{d}{dx}$. Thus, it contributes a term to the left and right of the Pascal Triangle in each differentiation step. For example, to obtain the second row of the operator $\delta_{4,3}$, we pass though the steps in Table 5, in that order, after each differentiation step.
\begin{align*} 
&\hspace{1.0cm}\ddx{4}{3} \\
&4\ddx{3}{3}+\ddx{4}{4}\\
12\ddx{2}{3}+&4\ddx{3}{4}+4\ddx{3}{4}+\ddx{4}{5}\\
24\xdx{3}+12\ddx{2}{4}+12\ddx{2}{4}+&4\ddx{3}{5}+12\ddx{2}{4}+4\ddx{3}{5}+4\ddx{3}{5}+\ddx{4}{6}
\end{align*}
\hspace*{5.8cm}{Table 5: The second row of the $\delta_{4,3}$\\                  

Table 6 corresponds to the finite sub-Pascal triangle of order four, $(r+1)$, with rows given in Table 5.

\hspace*{7.65cm} 1\\ 
\hspace*{7.75cm} 1 \hspace{.5cm} 1  \\
\hspace*{7.4cm} 1\hspace{.5cm} 2 \hspace{.5cm}1\\
\hspace*{7.2cm}1\hspace{.5cm}3\hspace{.6cm}3\hspace{.5cm}1\\\\
\hspace*{5.8cm}{Table 6: sub-Pascal triangle for $r=3$\\

\normalsize 

In the case of $\delta_{3,4}$, there will be no `two-pair' contributions after the third step. Thus, $\min\{r,m\}$ decides the length of the sub-Pascal triangle. This observation extends Theorem \ref{rec1} to the more general version given in the next result. 
\begin{theorem}\label{rec2}
Let $r,m\in \mathbb{N}=\{1, 2, 3, \cdots\}$ and $r\leq m$. The generalized Stirling numbers generated by the operators $\delta_{r,m}$ and $\delta_{m, r}$ are coincide with each other and follow the recurrence relation given by
\begin{equation}\label{eq-rec2}
S(n,k)=\sum_{i=0}^r \big(m+(m-r)(n-2)+k-i-1\big)_{r-i}\binom{r}{i} S(n-1,k-i), \;\;  0< k\leq n,
\end{equation}
with $S(0,0)=1$ and $S(n,0)=S(n,k)=0$ for $n>0$ and $1+min\{r,m\}(n-1) < k $ or $k\leq 0$. Here $(\cdot)_\lambda$ is the Pochhammer symbol for falling factorials. 
\end{theorem} 
\begin{proof}[Sketch of proof] 
The proof is similar to that of Theorem \ref{rec1} except that, in this case, it is required to consider $r+1$ many terms, which contribute to the coefficient $S(n,k)$. The appearance of the binomial coefficient is clear from the discussion of Section 4.1. The term $(m+(m-r)(n-2)+k-i-1)_{r-i}$ is obtained after $r-i$ many differentiations of the power of the basic term $x^{m+(m-r)(n-2)+k-i-1}\frac{d^{m+k-i}}{dx^{m+k-i}}$. This completes the proof.
\end{proof}

It is clear from Equation (\ref{eq-rec2}) that when $r=m=1$, the recurrence relation coincide with that of \emph{Stirling numbers of $2^{nd}$ kind} and when $r=2, m=1$, it coincide with that of the \emph{unsigned Lah numbers}. The Maple\textsuperscript{\small{\textregistered}} implementation of the recurrence relation (\ref{eq-rec2}) is given in the Algorithm \ref{gS}. 


All the sequences for the cases $k \in \{1,\, 2,\, 3, \cdots , 20\}$ are listed at A008277, A019538, A035342, A035469, A049029, A049385, A092082, A132056, A223511-A223522 \cite{key-21}. 

\section{New Class of Sequences}\label{seq}
We now consider \emph{rational} values of $r \in \mathbb{R}$. To begin our work, we first consider $r=\frac{1}{2}, \frac{1}{3}, \frac{1}{4}, \frac{1}{5}$ and $\frac{1}{6}$. 
\begin{definition}\label{dsthalf}
The generalized Stirling numbers of order $\frac{1}{2}$ is defined by the coefficients of $2^{\frac{n-1}{2}}\left(\sqrt{x}\frac{d}{dx}\right)^n$, when $n$ is odd-positive, and of $2^{\frac{n}{2}}\left(\sqrt{x}\frac{d}{dx}\right)^n$, when $n$ is even-positive. 
\end{definition}
The sequence corresponding to the Definition \ref{dsthalf} is given below, Sloane's A223168 \cite{key-21}.
\begin{table}[h]
\centering
\begin{tabular}{c|ccccccc}
\small{\backslashbox{n}{k}} &$c_{n,1}$ & $c_{n,2}$ & $c_{n,3}$ & $c_{n,4}$ & $c_{n,5}$ & $c_{n,6}$\\\hline
1 & 1 & &&&&\\
2 & 1 & 2 \\
3 & 3 & 2 \\
4 & 3 & 12 & 4\\
5 & 15 & 20 & 4 \\
6 & 15 & 90 & 60 & 8 \\
7 & 105 & 210 & 84 & 8\\
8 & 105 & 840 & 840 & 224 & 16\\
9 & 945 & 2520 & 1512 & 288 & 16\\
$\vdots$ & &&$\vdots$ & & & $\ddots$
\end{tabular}
\vspace{0.2in}
\begin{center}
Table 7. The sequence for $\delta_{\frac{1}{2},1}$ 
\end{center}
\end{table}

The expansion of this operator takes the form,
\begin{align*} 
2^0 \Big(\sqrt{x}\frac{d}{dx}\Big)^1 &= 1x^\frac{1}{2}\frac{d}{dx},\\
2^1 \Big(\sqrt{x}\frac{d}{dx}\Big)^2 &= 1\frac{d}{dx} + 2\xdx{2},\\
2^1 \Big(\sqrt{x}\frac{d}{dx}\Big)^3 &= 3x^\frac{1}{2}\frac{d^2}{dx^2} + 2x^\frac{3}{2}\frac{d^3}{dx^3},\\
2^2 \Big(\sqrt{x}\frac{d}{dx}\Big)^4 &= 3\frac{d^2}{dx^2} + 12\xdx{3} + 4\ddx{2}{4},\\
2^2 \Big(\sqrt{x}\frac{d}{dx}\Big)^5 &= 15x^\frac{1}{2}\frac{d^3}{dx^3} + 20x^\frac{3}{2}\frac{d^4}{dx^4} + 4x^\frac{5}{2}\frac{d^5}{dx^5},\\
&\vdots 
\end{align*}
\begin{center}
Table 8. Pseudo polynomials for the operator $\delta_{\frac{1}{2},1}$ 
\end{center}


These pseudo polynomials may be generated by the Maple code listed in Algorithm \ref{alg1}. According to V. Kotesovec \cite[A223168]{key-21}, the coefficients of the $k^{th}$ row of $\delta_{\frac{1}{2},1}$ is given by the $k^{th}$ derivative of $\frac{1}{2^{\lceil k/2 \rceil}}e^{x^2}$, where $\lceil \cdot \rceil$ is the ceiling function. The odd and even rows of this triangular sequence are listed at A223523 and A223524 \cite{key-21}. Odd rows have further interesting properties. For example, it includes absolute values of A098503 \cite{key-21} from right to left, which are coefficients of the \emph{generalized Laguerre polynomials} up to an integer multiple of $2^n$ for some $n\in \mathbb{N}$. They are given by the Rodrigues formula \cite{key-AS},

\begin{equation*}
L_n^{(\alpha)}(x)= {x^{-\alpha} e^x \over n!}{d^n \over dx^n} \left(e^{-x} x^{n+\alpha}\right),
\end{equation*}

which are solutions of the differential equation, 

\begin{equation*}
x\,y'' + (\alpha +1 - x)\,y' + n\,y = 0,
\end{equation*}

for arbitrary $\alpha$. The first few generalized Laguerre polynomials are given by,

\begin{align*}
L_0^{(\alpha)} (x) & = 1, \\
L_1^{(\alpha)}(x) & = -x + \alpha +1, \\
L_2^{(\alpha)}(x) & = \frac{x^2}{2} - (\alpha + 2)x + \frac{(\alpha+2)(\alpha+1)}{2}, \\
L_3^{(\alpha)}(x) & = \frac{-x^3}{6} + \frac{(\alpha+3)x^2}{2} - \frac{(\alpha+2)(\alpha+3)x}{2}
+ \frac{(\alpha+1)(\alpha+2)(\alpha+3)}{6},
\end{align*}

and the generalized Laguerre polynomial of degree $n$ is given by \cite{key-AS}

\begin{equation*}
L_n^{(\alpha)} (x) = \sum_{i=0}^n (-1)^i {n+\alpha \choose n-i} \frac{x^i}{i!}.
\end{equation*}
It can easily be seen that when $\alpha = 1/2$, we obtain A223523 \cite{key-21} up to signs of the terms. 

\subsection{Generalized $\delta$ Sequences}
The Definition \ref{dsthalf} can be extended for a large class of rational numbers. Let us now define $\delta_{r,1}$ sequences for other values of $r$ as well.
\begin{definition}\label{dstl}
The generalized Stirling numbers of order $r \in \{\frac{1}{2}, \frac{1}{3}, \frac{1}{4}, \ldots \}$, is defined by the coefficients of $r^{\frac{1-n}{2}}\left(x^r\frac{d}{dx}\right)^n$, when $n$ is odd-positive, and the coefficients of $r^{-\frac{n}{2}}\left(x^{1-r}\frac{d}{dx}\right)^n$, when $n$ is even-positive. 
\end{definition}

The corresponding sequence for $r \in \{\frac{1}{3}, \frac{1}{4}, \frac{1}{5}, \frac{1}{6}\}$ are listed in Table 9. All of these sequences and odd and even rows of them are listed at Sloane's A223168-A223172 and A223523-A223532 \cite{key-21}. 

\begin{remark}
The combinatorial interpretation of these sequences are unknown at the time of the production of this paper and the author suggests that they may be related to counting on Ferrer's-board, especially to \emph{r-creation boards}. The author's future work will be focused on the properties of these sequences and the representations of such sequences and beyond.
\end{remark}

\begin{remark}
It would be very interesting if one could classify all or a large class of such $\delta_{r,m}$ sequences. This could be a direction for future research. We would like to represent $S(n, k)$ as $\delta_{1,1}$-type and $L(n, k)$ as $\delta_{2,1}$-type, and etc. We may even consider \emph{negative integer} or \emph{irrational} values of $r \in \mathbb{R}$. To the author's knowledge this would lead to a new classification or at least a new method to represent some familiar sequences in combinatorial theory. 
\end{remark}

\section{Conclusion}


The paper presents the Mellin transforms of the generalized fractional integrals and derivatives, which generalize the Riemann-Liouville and the Hadamard fractional operators. 

We also show that the $\delta_{r,l}$ operator corresponds to several familiar sequences, specially \emph{generalized Stirling numbers} and \emph{Lah numbers}. We also show that the two operators $\delta_{r,m}$ and $\delta_{m,r}$, $r,m\in\mathbb{N}$, generate the same sequence. Further, we obtain the recurrence relations for the generalized Stirling numbers generated by the differential operators $\delta_{r,m}$ for $r,m\in\mathbb{N}$ and finally, we define a new class of sequences for $r \in \{\frac{1}{3}, \frac{1}{4}, \frac{1}{5}, \frac{1}{6}, \ldots \}$ and in tern show that $\delta_{\frac{1}{2},1}$ corresponds to the generalized Laguerre polynomials.

\begin{landscape}

\begin{table}[h]
	\centering
	  \subtable[$r= \frac{1}{3}$]{\begin{tabular}{c|ccccccc}
\small{\backslashbox{n}{k}} &$c_{n,1}$ & $c_{n,2}$ & $c_{n,3}$ & $c_{n,4}$ & $c_{n,5}$ & $c_{n,6}$\\\hline
1 & 1 & &&&&\\
2 & 1 & 3 \\
3 & 4 & 3 \\
4 & 4 & 24 & 9\\
5 & 28 & 42 & 9 \\
6 & 28 & 252 & 189 & 27 \\
7 & 280 & 630 & 270 & 27\\
8 & 280 & 3360 & 3780 & 1080 & 81\\
9 & 3640 & 10920 & 7020 & 1404 & 81\\
$\vdots$ & &&$\vdots$ & & & $\ddots$
\end{tabular}}
	  \hspace{0.5in}
	  \subtable[$r= \frac{1}{4}$]{\begin{tabular}{c|ccccccc}
\small{\backslashbox{n}{k}} &$c_{n,1}$ & $c_{n,2}$ & $c_{n,3}$ & $c_{n,4}$ & $c_{n,5}$ & $c_{n,6}$\\\hline
1 & 1 & &&&&\\
2 & 1 & 4 \\
3 & 5 & 4 \\
4 & 5 & 40 & 16\\
5 & 45 & 72 & 16 \\
6 & 45 & 540 & 432 & 64 \\
7 & 585 & 1404 & 624 & 64\\
8 & 585 & 9360 & 11232 & 3328 & 256\\
9 & 9945 & 31824 & 21216 & 4352 & 256\\
$\vdots$ & &&$\vdots$ & & & $\ddots$
\end{tabular}}\\
	  
\subtable[$r= \frac{1}{5}$]{\begin{tabular}{c|ccccccc}
\small{\backslashbox{n}{k}} &$c_{n,1}$ & $c_{n,2}$ & $c_{n,3}$ & $c_{n,4}$ & $c_{n,5}$ & $c_{n,6}$\\\hline
1 & 1 & &&&&\\
2 & 1 & 5 \\
3 & 6 & 5 \\
4 & 6 & 60 & 25\\
5 & 66 & 110 & 25 \\
6 & 66 & 990 & 825 & 125 \\
7 & 1056 & 2640 & 1200 & 125\\
8 & 1056 & 21120 & 26400 & 8000 & 625\\
9 & 22176 & 73920 & 50400 & 10500 & 625\\
$\vdots$ & &&$\vdots$ & & & $\ddots$
\end{tabular}}
	  \hspace{0.5in}
	  \subtable[$r= \frac{1}{6}$]{\begin{tabular}{c|ccccccc}
\small{\backslashbox{n}{k}} &$c_{n,1}$ & $c_{n,2}$ & $c_{n,3}$ & $c_{n,4}$ & $c_{n,5}$ & $c_{n,6}$\\\hline
1 & 1 & &&&&\\
2 & 1 & 6 \\
3 & 7 & 6 \\
4 & 7 & 84 & 36\\
5 & 91 & 156 & 36 \\
6 & 91 & 1638 & 1404 & 216 \\
7 & 1729 & 4446 & 2052 & 216\\
8 & 1729 & 41496 & 53352 & 16416 & 1296\\
9 & 43225 & 148200 & 102600 & 21600 & 1296\\
$\vdots$ & &&$\vdots$ & & & $\ddots$
\end{tabular}}
 \begin{center}
     Table 9. The generalized Stirling numbers for $r=\frac{1}{3}, \frac{1}{4}, \frac{1}{5}$, and $\frac{1}{6}.$   
\end{center}	  			
\end{table}
\end{landscape}



In a future project, we will derive formulae for the Laplace and Fourier transforms for the generalized fractional operators. We want to further investigate the effect on the new control parameter $\rho$, which makes more sense in this new setting.


\section*{Appendix}
This section contains the required Maple\textsuperscript{\small{\textregistered}}codes that generate the sequences we discuss throughout the paper. 
Some of these codes are also listed in the \emph{Sloane's On-line Encyclopedia of Integer Sequences} at \url{http://oeis.org} under the sequences A223168-A223172 and A223523-A223532 \cite{key-21}.

The pseudo polynomials of the operator $(x^r\frac{d^m}{dx^m})^n$ for $r,m\in \mathbb{N}$ are generated by the Maple code.
\begin{algorithm}\label{alg}
\begin{verbatim}
a[0] := f(x):
for i to 10 do 
   a[i] := simplify(x^r*diff(a[i-1],x$m)); 
end do
\end{verbatim}
\end{algorithm}

The generalized Stirling numbers $S(n,k)$ can be generated by the following Maple recursive procedure.

\begin{algorithm}\label{gS}
\begin{verbatim}
  genStrl := proc (n::integer, k::integer, r1::nonnegint, m1::nonnegint)
	    local r, m, i, total := 0; 
		    option remember;  
		    r := min(r1, m1); m := max(r1,m1); 
		     if n <= 0 or k <= 0 or 1+r*(n-1) < k then
			       return 0; 
			     elif n = 1 and k = 1 then
			       return 1
			     else 
			       for i from 0 to r do 
				         total := total+ mul(m+(m-r)*(n-2)+k-i-1-j, j = 0 .. r-i-1)
						                          * binomial(r, i)*genStrl(n-1, k-i, r, m);
				       end do; 
				       return total; 
			     end if; 
	end proc:
\end{verbatim}
\end{algorithm}

The pseudo polynomials for $r=\frac{1}{2}$ are generated by the Maple code given below. 
\begin{algorithm}\label{alg12}
\begin{verbatim}
   a[0] := f(x):                              
   for i from 1 to 10 do 
      a[i]:= simplify(2^((i+1)mod 2)*x^(1/2)*(diff(a[i-1],x$1)));
   end do;
\end{verbatim}
\end{algorithm}

The pseudo polynomials for $r=\frac{1}{\rho},\; \rho \in \{\frac{1}{3}, \frac{1}{4}, \frac{1}{5}, \ldots \}$  are generated by the Maple code in the Algorithm \ref{alg2}.
\begin{algorithm}\label{alg2}
\begin{verbatim}
   a[0] := f(x):                              
   for i from 1 to 10 do 
      a[i]:= simplify(r^((i+1)mod 2)*x^(((r-2)*((i+1)mod 2)+1)/r)
                                           *(diff(a[i-1],x$1)));
   end do;
\end{verbatim}
\end{algorithm}

The following recurrence relation based algorithm generates the generalized sequences in a triangular array. Here we give the case for $r=\frac{1}{4}$. Other cases are similar.
\newpage
\begin{algorithm}\label{alg1}\mbox{}
\begin{verbatim}
n := 10; r := 1/4; 
S := Array(1 .. n, 1 .. n); # This generates a nxn array of zeros
S[1, 1] := 1; 
 for i from 2 to n do 
   if i mod 2 = 1 then 
     for j to (i+1)*(1/2) do 
        S[i, j] := S[i-1,j]+j*S[i-1,j+1]  # i is odd
     end do 
   else 
     for j to (1/2)*i do 
        S[i, j] := i*(1+(j-1)/r)/(i-2*(j-1))*S[i-1,j] # i is even
     end do; 
        S[i,i/2+1] := 1/r*S[i-1,i/2] # For the case j = i/2+1
   end if 
 end do
\end{verbatim}
\end{algorithm}

\section*{Acknowledgment}
The author would like to thank Philip Feinsilver, the Department of Mathematics at Southern Illinois University, for pointing out interesting results related to Umbral Calculus, Darin B. Johnson, the Department of Mathematics at Delaware State University, for pointing out interesting connections to some literature in combinatorial theory and Jerzy Kocik, the Department of Mathematics at Southern Illinois University for valuable comments and suggestions. The author also would like to thank the anonymous referees for their invaluable suggestions, which put the article in its present form. 







\section*{References}



\end{document}